\documentclass[12pt]{report}
\usepackage{amssymb,latexsym,epsfig,mathrsfs,graphpap,graphics,amsmath,amssymb}
\usepackage{amstext}
\usepackage{amsthm}
\usepackage[utf8]{inputenc}
\usepackage[english]{babel}
\usepackage{helvet}
\usepackage{courier}

\theoremstyle{Definition}

\theoremstyle{remark}

\numberwithin{equation}{section}

\begin{document}

\begin{center}

 {\bf\Large {Inequalities for Nonuniform Wavelet Frames }}

\parindent=0mm \vspace{.4in}
{\bf{ Firdous A. Shah$^{\star}$}}

\end{center}

\begin{quote}
\parindent=0mm \vspace{.1in}
{\it\small{$^{\star}$Department of  Mathematics, University of Kashmir,  South Campus, Anantnag-192101,  Jammu and Kashmir, India. E-mail:
 $\text{fashah79@gmail.com}$}}

\parindent=0mm \vspace{.1in}
 { \small{\bf Abstract:} Gabardo and Nashed have studied  nonuniform wavelets based on the theory of spectral pairs  for which the associated translation set $\Lambda =\left\{ 0,r/N\right\}+2\,\mathbb Z$ is no longer a discrete subgroup of $\mathbb R$ but a spectrum associated with a certain one-dimensional spectral pair and the associated dilation is an even positive integer related to the given spectral pair. In this paper, we construct the associated wavelet frames and establish some  sufficient conditions that ensure the nonuniform wavelet system $\left\{ \psi_{j, \lambda}(x)= (2N)^{j/2}\psi\big( (2N)^{j}x-\lambda\big), j\in \mathbb Z, \lambda\in \Lambda \right\}$  to be a frame for  $L^2(\mathbb R)$. The conditions proposed are stated in terms of the Fourier transforms of the wavelet system's generating functions.}

\parindent=0mm \vspace{.1in}

\parindent=0mm \vspace{.1in}

\end{quote}

\parindent=0mm \vspace{.1in}
{\bf{1. Introduction}}

\parindent=0mm \vspace{.1in}
Frames were first introduced by Duffin and Schaeffer [6] to investigate some deep problems in nonharmonic Fourier series, and more particularly with the question of determining when a family of exponentials $\left\{e^{i\alpha_{n}t}:n\in\mathbb Z\right\}$ is complete for $L^2[a, b]$. Obviously, the significance of the idea was not understood by the mathematical community; at least it took nearly thirty years before the next treatment appeared in print. In 1986, Daubechies et al.[4] reintroduced frames and observed that frames can be used to find series expansions of function in $L^2(\mathbb R)$, which are similar to the expansion using orthonormal basis. Since then, the theory of frames began to be studied widely and deeply. The redundancy and flexibility offered by frames has impelled their application in several areas of mathematics, physics and engineering [1,5]. Mathematically, a sequence $\left\{f_{k}\right\}_{k=1}^{\infty}$ of elements  of a Hilbert space $\mathcal H$ is called a {\it frame} for $\mathcal H$ if there exist constants $A$ and $B$ with $0<A\le B<\infty$ such that
\begin{align*}
A\big\|f\big\|_{2}^{2}\le \sum_{k=1}^{\infty}\left|\big\langle f, f_{k}\big\rangle \right|^2\le B \big\|f\big\|_{2}^{2},\quad \text{ for all}~f\in \mathcal H.\tag{1.1}
\end{align*}

\parindent=0mm \vspace{.0in}
 The greatest possible such $A$ is the lower frame bound and the least possible such $B$ is the upper frame bound. A tight frame refers to the case when $A = B$, and a Parseval frame refers to the case when $A = B = 1$.

\parindent=8mm \vspace{.1in}
The particular frames of interest to us will be the frames in the space $L^2(\mathbb R)$ which are generated by the action of dilations and translations on a single function or finite number of functions in $L^2(\mathbb R)$. In order to describe these frames, we define  wavelet systems of the form
\begin{align*}
{\cal W}(\psi, a, b)=\Big\{\psi_{j,k}=: a^{j/2}\psi\big(a^{j}x-kb\big):j,k\in\mathbb Z\Big\},\tag{1.2}
\end{align*}

\parindent=0mm \vspace{.0in}
where $\psi\in L^2(\mathbb R)$ and  $a, b\in \mathbb R$ with $a>1,$ and $b>0$. Wavelet systems ${\cal W}(\psi, a, b)$  that constitute frames  for $L^2(\mathbb R)$ have a wide variety of applications.  Daubechies [3] proved the first result on the necessary and sufficient conditions for the wavelet system ${\cal W}(\psi, a, b)$ to be frame for $L^2(\mathbb R)$, and since then, these conditions have been improved and investigated by many authors [1,2,10,14].

\parindent=8mm \vspace{.2in}
 All these concepts are developed on regular lattices, that is the translation set is always a group. Recently, Gabardo and  Nashed [8,9] developed the theory of nonuniform wavelets and wavelet sets in $L^2(\mathbb R)$ for  which the translation set is no longer a discrete subgroup  of $\mathbb R$, but a union of two lattices, which is associated with a famous open conjecture of Fuglede on spectral pairs [7]. The main results of Gabardo and Nashed deal with necessary and sufficient condition for the existence of associated wavelets and extension of Cohen’s theorem.  Sharma and Manchanda [15] presented a necessary and sufficient conditions for nonuniform wavelet frames in the frequency domain. The theory of nonuniform wavelets was further studied and investigated by several authors in different directions including wavelets, vector-valued wavelets  and wavelet packets on local fields of positive characteristic [11--13].

\parindent=8mm \vspace{.2in}
One of the fundamental problems in the study of wavelet frames is to find conditions on the wavelet function and the dilation and translation parameters so that the corresponding wavelet system forms a frame for $L^2(\mathbb R)$. Our purpose in the paper is to introduce and establish conditions for the nonuniform wavelet systems associated with spectral pairs to be frames for $L^2(\mathbb R)$. More precisely, we derive three sufficient conditions under which nonuniform wavelet system $\left\{ \psi_{j, \lambda}(x)= (2N)^{j/2}\psi\big( (2N)^{j}x-\lambda\big), j\in \mathbb Z, \lambda\in \Lambda \right\}$  become frame in $L^2(\mathbb R)$. The inequalities we proposed are stated in terms of the Fourier transforms of the wavelet system's generating functions and all  these result are valid without any decay assumptions on the generators of the system.

\parindent=8mm \vspace{.2in}
The paper is structured as follows. In Section 2, we introduce some notations and preliminaries related to the nonuniform wavelets associated with
one-dimensional spectral pairs. Section 3 is devoted to the discussion of sufficient conditions for nonuniform wavelet frame  in $L^2(\mathbb R)$,
 and three such conditions are given.

\pagestyle{myheadings}

\parindent=0mm \vspace{.2in}

{\bf{2. Nonuniform Wavelet Frames in $L^2(\mathbb R)$ }}

\parindent=0mm \vspace{.2in}
For an integer $N\ge 1$ and an odd integer $r$ with $1 \le r \le 2N-1$ such that $r$ and $N$ are relatively prime, we define

$$\Lambda=\left\{0, \dfrac{r}{N}\right\}+2\mathbb Z=\left\{ \dfrac{rk}{N}+2n: n\in \mathbb Z, k=0,1\right\}.\eqno(2.1)$$

\parindent=0mm \vspace{.1in}
It is easy to verify that $\Lambda$ is not necessarily a group nor a uniform discrete set, but is the union of $\mathbb Z$ and a translate of $\mathbb Z$.
 Moreover, the set $\Lambda$ is the spectrum for the spectral set $\Upsilon=\big[0, \frac{1}{2}\big)\cup \big[\frac{N}{2}, \frac{N+1}{2}\big)$ and the
pair $(\Lambda,\Upsilon)$ is called a {\it spectral pair} [7,8].

\parindent=8mm \vspace{.2in}
For a given   $\psi \in L^2(\mathbb R)$, define the nonuniform wavelet system
\begin{align*}
{\cal N}(\psi, j, \lambda)=\Big\{\psi_{j,\lambda}=: (2N)^{j/2}\psi\big((2N)^{j}x-\lambda\big):j\in\mathbb Z, \lambda\in \Lambda\Big\}.\tag{2.2}
\end{align*}

\parindent=0mm \vspace{.0in}
We call the wavelet system ${\cal N}(\psi, j, \lambda)$	 a nonuniform wavelet frame for $L^2(\mathbb R)$, if there exist positive numbers
$0 < C \le D < \infty$ such that
\begin{align*}
C\big\|f\big\|_{2}^{2}\le \sum_{j\in\mathbb Z}\sum_{\lambda\in \Lambda} \left|\big\langle f, \psi_{j,\lambda}\big\rangle\right|^2 \le D\big\|f\big\|_{2}^{2},\tag{2.3}
\end{align*}

\parindent=0mm \vspace{.0in}
holds for every $f\in  L^2(\mathbb R)$. In order to prove our main results, we need the following lemma whose proof can be found in [1].

\parindent=0mm \vspace{.1in}
{\bf{ Lemma 2.1.}} Suppose that $\left\{f_{k}\right\}_{k=1}^{\infty}$ is a family of elements in a Hilbert space $\mathcal H$ such that the inequalities
(1.1) holds for all $f$ in a dense subset ${\cal D} $ of $\mathcal H$.
Then, the same inequalities is true for all $f\in \mathcal H$.

\parindent=8mm \vspace{.1in}
In view of Lemma 2.1, we will consider the following set of functions:
\begin{align*}
{\cal D}= \left\{ f\in L^2(\mathbb R):\hat f\in L^{\infty}(\mathbb R) ~\text{and}~\hat f~\text{has compact support in}~\mathbb R\right\}.
\end{align*}

\parindent=0mm \vspace{.0in}
It is clear that ${\cal D}$ is a dense subspace of $L^2(\mathbb R)$. Therefore, it is enough to verify that the system ${\cal N}(\psi, j, \lambda)$	 given by (2.2) is a frame for $L^2(\mathbb R)$ if (2.3)  hold for all $f\in {\cal D}$. Moreover, we need the following lemma on nonuniform wavelet frames whose proof  can be found in [15, Lemma 3.1].

\parindent=0mm \vspace{.1in}
{\bf{ Lemma 2.2.}} {\it Let $f\in {\cal D}$ and $\psi \in L^2(\mathbb R)$. If esssup$\big\{\sum_{j \in\mathbb Z}\big|\hat\psi \big(
(2N)^{-j}\xi \big)\big|^2 : \xi\in[1, 2N]\big\}<\infty$, then }
$$\sum_{j\in\mathbb Z}\sum_{\lambda\in \Lambda} \left|\big\langle f, \psi_{j,\lambda}\big\rangle\right|^2 =\int_{\mathbb R}\left|\hat f(\xi)\right|^{2}
\sum_{j \in\mathbb Z}\left|\hat\psi \big((2N)^{-j}\xi \big)\right|^2 d\xi+R_{\psi}(f),\eqno(2.4)$$

\parindent=0mm \vspace{.1in}
{\it where $R_{\psi}(f)=R_{0}+R_{1}+\dots+R_{2N-1}$ and for $0 \le p \le 2N-1, R_{p}$ is given by}
\begin{align*}
R_{p}&=\dfrac{1}{4N}\displaystyle\sum_{j \in \mathbb Z}\sum_{\ell\ne p}\int_{\mathbb R}
\left\{\overline{\hat f\left(\xi+(2N)^j\dfrac{p}{2}\right)}\hat\psi\left(\dfrac{\xi}{(2N)^j}
+\dfrac{p}{2}\right)\hat f\left(\xi+(2N)^j\dfrac{\ell}{2}\right)\right.\\
&\qquad\qquad\qquad\qquad\qquad\qquad\qquad\times\left.\overline{\hat\psi\left(\dfrac{\xi}{(2N)^j}+\dfrac{\ell}{2}\right)}
\left(1+e^{\pi i\frac{r}{N}(\ell-p)}\right)\right\}d\xi.\tag{2.5}
\end{align*}

\parindent=8mm \vspace{.0in}
We now establish our first sufficient condition for the  nonuniform wavelet system ${\cal N}(\psi, j, \lambda)$ given by (2.2) to be a frame
for $L^2(\mathbb R)$. For this, we set
\begin{align*}
\Delta_\psi(m)=\text{ess}\sup_{}\left\{\sum_{j \in \mathbb Z}\left|t_\psi\left(m, \dfrac{\xi}{(2N)^j}\right)\right|:\xi\in[1, 2N]\right\},\tag{2.6}
\end{align*}

\parindent=0mm \vspace{.0in}
where
\begin{align*}
t_\psi\left(m, \xi\right)=\sum_{k \in \mathbb N_0}\hat\psi\Big((2N)^k\xi\Big)\overline{\hat\psi\left((2N)^k\Big(\xi+\dfrac{m}{2}\Big)\right)}.\tag{2.7}
\end{align*}

\parindent=0mm \vspace{.0in}
We also use the following set:

$$\Omega=\Big\{ (2N)k+\ell: k\in\mathbb N_{0}, 1\le \ell\le 2N-1 \Big\}.$$

\parindent=0mm \vspace{.1in}
Analogous to the uniform case, we give the first sufficient condition as follows.

\parindent=0mm \vspace{.2in}
{\bf Theorem 2.3.} {\it Suppose $\psi \in L^2(\mathbb R)$ such that}
\begin{align*}
A_\psi&=\text{ess}\inf_{\xi\in[1, 2N]}\sum_{j \in\mathbb Z}\left|\hat\psi\left(\dfrac{\xi}{(2N)^j}\right)\right|^2-\sum_{\alpha \ne \beta
 \in \Omega}\left[\Delta_\psi\left(\dfrac{\alpha - \beta}{2}\right)\cdot\Delta_\psi\left(\dfrac{\beta-\alpha}{2}\right)\right]^{1/2}>0,\\
B_\psi&=\text{ess}\sup_{\xi\in[1, 2N]}\sum_{j \in\mathbb Z}\left|\hat\psi\left(\dfrac{\xi}{(2N)^j}\right)\right|^2+
\sum_{\alpha \ne \beta \in \Omega}\left[\Delta_\psi\left(\dfrac{\alpha - \beta}{2}\right)\cdot\Delta_\psi\left(\dfrac{\beta-\alpha}{2}\right)\right]^{1/2}
<\infty.
\end{align*}

\parindent=0mm \vspace{.0in}
{\it Then $\big\{\psi_{j,\lambda}:j \in \mathbb Z,\lambda \in \Lambda\big\}$ is a frame for $L^2(\mathbb R)$ with bounds $A_\psi$ and $B_\psi$.}

\parindent=0mm \vspace{.2in}
{\it Proof.}  Since the last series in (2.5) is absolutely convergent for every $f\in {\cal D}$, we can estimate $R_{\psi}(f)$ by
rearranging the series, changing the orders of summation and integration by Levi Lemma so that we deduce that
\begin{align*}
\Big|R_\psi(f)\Big|&\le\dfrac{1}{2N}\displaystyle\sum_{p=0}^{2N-1}\sum_{j \in \mathbb Z}\int_{\mathbb R}\left|\overline{\hat f(\xi)}\hat\psi
\left(\dfrac{\xi}{(2N)^j}\right)\right|\left\{\sum_{\ell\ne p}\left|\hat f\left(\xi+(2N)^j\left(\dfrac{\ell-p}{2}\right)\right)\right.\right.\
&\qquad\qquad\qquad\qquad\qquad\qquad\qquad\qquad\qquad\qquad\left.\left.\times\overline{\hat\psi\left(\dfrac{\xi}{(2N)^j}+
\dfrac{\ell-p}{2}\right)}\right|\right\}d\xi\\
&=\dfrac{1}{2N}\displaystyle\sum_{\alpha=0}^{2N-1}\sum_{j \in \mathbb Z}\int_{\mathbb R}\left|\overline{\hat f(\xi)}\right|
\left\{\sum_{k\in \mathbb N_0}\sum_{\alpha \ne \beta\in \Omega}\left|\hat\psi\left(\dfrac{\xi}{(2N)^j}\right)\right.\right.\\
&\qquad\qquad\left.\left.\times\hat f\left(\xi+(2N)^{j+k}\left(\dfrac{\alpha - \beta}{2}\right)\right)\overline{\hat\psi
\left(\dfrac{\xi}{(2N)^j}+(2N)^k\left(\dfrac{\alpha - \beta}{2}\right)\right)}\right|\right\}d\xi\\
&=\dfrac{1}{2N}\displaystyle\sum_{\alpha=0}^{2N-1}\int_{\mathbb R}\left|\overline{\hat f(\xi)}\right|\left\{\sum_{k\in \mathbb N_0}
\sum_{\alpha \ne \beta\in \Omega}\sum_{j \in \mathbb Z}\left|\hat\psi\left(\dfrac{\xi}{(2N)^{j-k}}\right)\right.\right.\\
&\qquad\qquad\left.\left.\times\hat f\left(\xi+(2N)^{j}\left(\dfrac{\alpha - \beta}{2}\right)\right)\overline{\hat\psi\left(\dfrac{\xi}{(2N)^{j-k}}+(2N)^k
\left(\dfrac{\alpha - \beta}{2}\right)\right)}\right|\right\}d\xi\\
&=\dfrac{1}{2N}\displaystyle\sum_{\alpha=0}^{2N-1}\int_{\mathbb R}\left|\overline{\hat f(\xi)}\right|\left\{\sum_{j \in \mathbb Z}
\sum_{\alpha \ne \beta\in \Omega}\left|\hat f\left(\xi+(2N)^{j}\left(\dfrac{\alpha - \beta}{2}\right)\right)\right.\right.\\
&\qquad\qquad\qquad\qquad\left.\left.\times\displaystyle\sum_{k\in \mathbb N_{0}}\hat\psi\left(\dfrac{\xi}{(2N)^{j-k}}\right)
\overline{\hat\psi\left((2N)^{k}\left(\dfrac{\xi}{(2N)^j}+\dfrac{\alpha - \beta}{2}\right)\right)}\right|\right\}d\xi\\
&=\dfrac{1}{2N}\displaystyle\sum_{\alpha=0}^{2N-1}\int_{\mathbb R}\left|\overline{\hat f(\xi)}\right|\left\{\sum_{j \in \mathbb Z}
\sum_{\alpha \ne \beta\in \Omega}\left|\hat f\left(\xi+(2N)^{j}\left(\dfrac{\alpha - \beta}{2}\right)\right)\right.\right.\\
&\qquad\qquad\qquad\qquad\qquad\qquad\qquad\qquad\qquad\qquad\left.\times\left| t_\psi\left(\dfrac{\alpha - \beta}{2}, \dfrac{\xi}{(2N)^j}\right)
\right|\right\}d\xi\\
&=\dfrac{1}{2N}\displaystyle\sum_{\alpha=0}^{2N-1}\sum_{j \in \mathbb Z}\sum_{\alpha \ne \beta\in \Omega}\int_{\mathbb R}
\left\{\left|\hat f(\xi)\right|\left| t_\psi\left(\dfrac{\alpha - \beta}{2}, \dfrac{\xi}{(2N)^j}\right)\right|^{1/2}\right\}
\end{align*}
\begin{align*}
&\qquad\qquad\qquad\qquad\times\left\{\left|\hat f\left(\xi+(2N)^{j}\left(\dfrac{\alpha - \beta}{2}\right)\right)\right|
\left| t_\psi\left(\dfrac{\alpha - \beta}{2}, \dfrac{\xi}{(2N)^j}\right)\right|^{1/2}\right\}d\xi\\
&\le\dfrac{1}{2N}\displaystyle\sum_{\alpha=0}^{2N-1}\sum_{j \in \mathbb Z}\sum_{\alpha \ne \beta\in \Omega}\left\{\int_{\mathbb R}\left|\hat f(\xi)\right|^2
\left| t_\psi\left(\dfrac{\alpha - \beta}{2}, \dfrac{\xi}{(2N)^j}\right)\right|d\xi\right\}^{1/2}\\
&\qquad\qquad\qquad\times\left\{\displaystyle\int_{\mathbb R}\left|\hat f\left(\xi+(2N)^{j}\left(\dfrac{\alpha - \beta}{2}\right)\right)\right|^2
\left| t_\psi\left(\dfrac{\alpha - \beta}{2}, \dfrac{\xi}{(2N)^j}\right)\right|d\xi\right\}^{1/2}\\
&\le\dfrac{1}{2N}\displaystyle\sum_{\alpha=0}^{2N-1}\sum_{\alpha \ne \beta\in \Omega}\left\{\sum_{j \in \mathbb Z}
\int_{\mathbb R}\left|\hat f(\xi)\right|^2\left| t_\psi\left(\dfrac{\alpha - \beta}{2}, \dfrac{\xi}{(2N)^j}\right)\right|d\xi\right\}^{1/2}\\
&\qquad\qquad\qquad\times\left\{\displaystyle\sum_{j \in \mathbb Z}\int_{\mathbb R}\left|\hat f\left(\xi+(2N)^{j}\left(\dfrac{\alpha - \beta}{2}\right)\right)
\right|^2\left| t_\psi\left(\dfrac{\alpha - \beta}{2}, \dfrac{\xi}{(2N)^j}\right)\right|d\xi\right\}^{1/2}\\
&=\dfrac{1}{2N}\displaystyle\sum_{\alpha=0}^{2N-1}\sum_{\alpha \ne \beta\in \Omega}\left\{\sum_{j \in \mathbb Z}\int_{\mathbb R}
\left|\hat f(\xi)\right|^2\left| t_\psi\left(\dfrac{\alpha - \beta}{2}, \dfrac{\xi}{(2N)^j}\right)\right|d\xi\right\}^{1/2}\\
&\qquad\qquad\qquad\qquad\qquad\qquad\qquad\times\left\{\displaystyle\sum_{j \in \mathbb Z}\int_{\mathbb R}\left|\hat f(\omega)\right|^2
\left| t_\psi\left(-\dfrac{\alpha - \beta}{2}, \dfrac{\omega}{(2N)^j}\right)\right|d\xi\right\}^{1/2}\\
&\le\dfrac{1}{2N}\displaystyle\sum_{\alpha=0}^{2N-1}\sum_{\alpha \ne \beta\in \Omega}\left\{\int_{\mathbb R}\left|\hat f(\xi)\right|^2
\Delta_\psi\left(\dfrac{\alpha - \beta}{2}\right)d\xi\right\}^{1/2}\\
&\qquad\qquad\qquad\qquad\qquad\qquad\qquad\qquad\qquad\times\left\{\displaystyle\int_{\mathbb R}\left|\hat f(\xi)\right|^2
\Delta_\psi\left(\dfrac{-(\alpha - \beta)}{2}\right)d\xi\right\}^{1/2}\\
&=\dfrac{1}{2N}\displaystyle\sum_{\alpha=0}^{2N-1}\int_{\mathbb R}\left|\hat f(\xi)\right|^2 d\xi\sum_{\alpha \ne \beta\in \Omega}
\left[\Delta_\psi\left(\dfrac{\alpha - \beta}{2}\right)\cdot \Delta_\psi\left(\dfrac{-(\alpha - \beta)}{2}\right)\right]^{1/2}.
\end{align*}

\parindent=0mm \vspace{.1in}
Consequently, it follows from the expression (2.4) in Lemma 2.2 that
\begin{align*}
\displaystyle\sum_{j \in \mathbb Z}\sum_{\lambda \in \Lambda}\left|\big\langle f, \psi_{j, \lambda}\big\rangle\right|^2&\ge
 \displaystyle\int_{\mathbb R}\left|\hat f(\xi)\right|^2\left\{\sum_{j \in \mathbb Z}\left|\hat \psi\left(\dfrac{\xi}{(2N)^j}\right)\right|^2\right.\\\
&\quad\left.-\displaystyle\sum_{\alpha \ne \beta\in \Omega} \left[\Delta_\psi\left(\dfrac{\alpha - \beta}{2}\right)\cdot \Delta_\psi
\left(\dfrac{-(\alpha - \beta)}{2}\right)\right]^{1/2}\right\}d\xi,\tag{2.8}
\end{align*}

\parindent=0mm \vspace{.0in}
and
\begin{align*}
\displaystyle\sum_{j \in \mathbb Z}\sum_{\lambda \in \Lambda}\left|\big\langle f, \psi_{j, \lambda}\big\rangle\right|^2&\le \displaystyle\int_{\mathbb R}
\left|\hat f(\xi)\right|^2\left\{\sum_{j \in \mathbb Z}\left|\hat \psi\left(\dfrac{\xi}{(2N)^j}\right)\right|^2\right.\\\
&\quad\left.+\displaystyle\sum_{\alpha \ne \beta\in \Omega} \left[\Delta_\psi\left(\dfrac{\alpha - \beta}{2}\right)\cdot \Delta_\psi
\left(\dfrac{-(\alpha - \beta)}{2}\right)\right]^{1/2}\right\}d\xi.\tag{2.9}
\end{align*}

\parindent=0mm \vspace{.1in}
Taking infimum in (2.8) and supremum in (2.9), respectively, we obtain that
\begin{align*}
A_\psi\big\|f\big\|^2_2\le \sum_{j \in \mathbb Z}\sum_{\lambda \in \Lambda}\left|\big\langle f, \psi_{j, \lambda}\big\rangle\right|^2\le B_\psi\big\|f\big\|^2_2,
\end{align*}

\parindent=0mm \vspace{.0in}
hold for all $f\in {\cal D}$. This completes the proof of Theorem 2.3.\qquad\fbox

\parindent=8mm \vspace{.2in}

Before, we state our next sufficient condition, we introduce some notations. Similar to the $a$-adic number, we call an element $\alpha \in \mathbb R$,  a $2N$-adic
number if it has the form $\alpha=(2N)^{j}(\lambda-\sigma)/2,\,j \in \mathbb Z,\,\lambda \ne \sigma \in \Lambda$. With this concept, we consider the set
$$\Gamma=\left\{\alpha \in \mathbb R: {\text{there exists}}\;(j,\lambda)\in \mathbb Z\times \Lambda\;{\text{such that}}\;\alpha=\dfrac{(2N)^{j}(\lambda-\sigma)}{2};
\lambda\ne\sigma\right\},\eqno(2.10)$$

\parindent=0mm \vspace{.1in}
and for all $\alpha \in \Gamma$, we define
\begin{align*}
\displaystyle I(\alpha)&=\left\{(j,\lambda)\in \mathbb Z\times \Lambda:\alpha =\dfrac{(2N)^{j}(\lambda-\sigma)}{2}
\right\},\tag{2.11}\\\
\displaystyle \Delta_\alpha^+(\xi) &=\sum_{(j,\lambda)\ne(j,\sigma)\in I(\alpha)}\hat\psi\left(\dfrac{\xi}{(2N)^j}\right)\overline{\hat\psi\left(\dfrac{\xi}{(2N)^j}
+\dfrac{\lambda-\sigma}{2}\right)},\tag{2.12}\\\
\Delta_\alpha^-(\xi) &=\sum_{(j,\lambda)\ne(j,\sigma)\in I(\alpha)}\hat\psi\left(\dfrac{\xi}{(2N)^j}\right)\overline{\hat\psi\left(\dfrac{\xi}{(2N)^j}
-\dfrac{\lambda-\sigma}{2}\right)}.\tag{2.13}
\end{align*}

\parindent=0mm \vspace{.1in}
With the notations above we state the following result.

\parindent=0mm \vspace{.2in}
{\bf {Theorem 2.4.}} {\it Suppose $\psi \in L^2(\mathbb R)$ such that}
\begin{align*}
C_\psi&=ess\inf_{\xi\in[1, 2N]}\Big\{\Delta_0^+(\xi)-\sum_{\alpha \in \Gamma\setminus\{0\}}\left|\Delta_\alpha^+(\xi)\right|\Big\}>0,\\
D_\psi&=ess\sup_{\xi\in[1, 2N]}\sum_{\alpha \in \Gamma\setminus\{0\}}\left|\Delta_\alpha^+(\xi)\right|<+\infty.
\end{align*}

\parindent=0mm \vspace{.0in}
{\it Then $\big\{\psi_{j,\lambda}:j \in \mathbb Z,\,\lambda \in \Lambda\big\}$ is a wavelet frame for $L^2(\mathbb R)$ with bounds $C_\psi$ and $D_\psi$.}

\parindent=0mm \vspace{.2in}
{\it Proof.} We first note that  $\Delta_0^+(\xi)=\sum_{j \in\mathbb Z}\big|\hat\psi\left(\xi/(2N)^j\right)\big|^2$ by the definition of $\Delta_\alpha^+(\xi)$.
 We apply Lemma 2.2 to re-estimate $R_\psi(f)$ for $f\in {\cal D}$ as
\begin{align*}
\big|R_\psi(f)\big|&=\left|\dfrac{1}{4N}\displaystyle\sum_{p=0}^{2N-1}\sum_{j \in \mathbb Z}\sum_{\ell\ne p}\int_{\mathbb R}
\left\{\overline{\hat f\left(\xi+(2N)^j\dfrac{p}{2}\right)}\hat\psi\left(\dfrac{\xi}{(2N)^j}
+\dfrac{p}{2}\right)\hat f\left(\xi+(2N)^j\dfrac{\ell}{2}\right)\right.\right.\\
&\qquad\qquad\qquad\qquad\qquad\qquad\qquad\left.\times\left.\overline{\hat\psi\left(\dfrac{\xi}{(2N)^j}+\dfrac{\ell}{2}\right)}
\left(1+e^{\pi i\frac{r}{N}(\ell-p)}\right)\right\}d\xi\right|\\
&\le\dfrac{1}{2N}\displaystyle\sum_{p=0}^{2N-1}\sum_{j \in \mathbb Z}\sum_{\ell\ne p}\int_{\mathbb R}\left|\overline{\hat f\left(\xi+(2N)^j\dfrac{p}{2}\right)}\hat f\left(\xi+(2N)^j\dfrac{\ell}{2}\right)\right|
\end{align*}
\begin{align*}
&\qquad\qquad\qquad\qquad\qquad\qquad\qquad\times\left|\hat\psi\left(\dfrac{\xi}{(2N)^j}+\dfrac{p}{2}\right)\overline{\hat\psi\left(\dfrac{\xi}{(2N)^j}
+\dfrac{\ell}{2}\right)}\right|d\xi\\
&=\displaystyle\sum_{\alpha \in \Gamma\setminus\{0\}}\sum_{(j, \lambda)\ne (j, \sigma)\in I(\alpha)}\int_{\mathbb R}\left|
\overline{\hat f\left(\xi+(2N)^j\dfrac{\lambda}{2}\right)}\hat f\left(\xi+(2N)^j\dfrac{\sigma}{2}\right)\right|\\\
&\qquad\qquad\qquad\qquad\qquad\qquad\qquad\times\left|\hat\psi\left(\dfrac{\xi}{(2N)^j}+\dfrac{\lambda}{2}\right)\overline{\hat\psi\left(\dfrac{\xi}{(2N)^j}
+\dfrac{\sigma}{2}\right)}\right|d\xi\\
&\le\displaystyle\sum_{\alpha \in \Gamma\setminus\{0\}}\sum_{(j, \lambda)\ne (j, \sigma)\in I(\alpha)}\int_{\mathbb R}\left|\overline{\hat f(\xi)}
\hat f\left(\xi+(2N)^j\left(\dfrac{\lambda-\sigma}{2}\right)\right)\right|\\
&\qquad\qquad\qquad\qquad\qquad\qquad\qquad\times\left|\hat\psi\left(\dfrac{\xi}{(2N)^j}\right)\overline{\hat\psi\left(\dfrac{\xi}{(2N)^j}
+\dfrac{\lambda-\sigma}{2}\right)}\right|d\xi\\
&=\displaystyle\sum_{\alpha \in \Gamma\setminus\{0\}}\int_{\mathbb R}\left|\overline{\hat f(\xi)}\hat f\left(\xi+\dfrac{\alpha}{2}\right)\right|
\left\{\displaystyle\sum_{(j, \lambda)\ne (j, \sigma)\in I(\alpha)}\left|\hat\psi\left(\dfrac{\xi}{(2N)^j}\right)\overline{\hat\psi
\left(\dfrac{\xi}{(2N)^j}+\dfrac{\lambda-\sigma}{2}\right)}\right|\right\}d\xi\\
&=\displaystyle\sum_{\alpha \in \Gamma\setminus\{0\}}\int_{\mathbb R}\left|\overline{\hat f(\xi)}\hat f\left(\xi+\dfrac{\alpha}{2}\right)\right|
\Big|\Delta_\alpha^+(\xi)\Big|d\xi\qquad {\big(\text{By Eq. (2.12)}\big)}\\
&\le\displaystyle\sum_{\alpha \in \Gamma\setminus\{0\}}\left\{\int_{\mathbb R}\left|\hat f(\xi)\right|^2\left|\Delta_\alpha^+(\xi)\right|d\xi\right\}^{1/2}
\left\{\int_{\mathbb R}\left|\hat f\left(\xi+\dfrac{\alpha}{2}\right)\right|^2\left|\Delta_\alpha^+(\xi)\right|d\xi\right\}^{1/2}\\
&\le\left\{\displaystyle\sum_{\alpha \in \Gamma\setminus\{0\}}\int_{\mathbb R}\left|\hat f(\xi)\right|^2\left|\Delta_\alpha^+(\xi)\right|d\xi\right\}^{1/2}
\left\{\displaystyle\sum_{\alpha \in \Gamma\setminus\{0\}}\int_{\mathbb R}\left|\hat f\left(\xi+\dfrac{\alpha}{2}\right)\right|^2
\left|\Delta_\alpha^+(\xi)\right|d\xi\right\}^{1/2}.\tag{2.14}
\end{align*}

\parindent=0mm \vspace{.1in}
Let $\omega=\xi+{\alpha}/{2}$. We deduce from  $\alpha=(2N)^{j}(\lambda-\sigma)$ for $(j, \lambda)\ne (j,\sigma)\in I(\alpha)$ that
\begin{align*}
\Delta_\alpha^+(\xi)&=\displaystyle\sum_{(j,\lambda)\ne(j,\sigma)\in I(\alpha)}\hat\psi\left(\dfrac{\xi}{(2N)^j}\right)
\overline{\hat\psi\left(\dfrac{\xi}{(2N)^j}+\dfrac{\lambda-\sigma}{2}\right)}\\
&=\displaystyle\sum_{(j,\lambda)\ne(j,\sigma)\in I(\alpha)}\hat\psi\left((2N)^{-j}\left(\omega-\dfrac{\alpha}{2}\right)\right)
\overline{\hat\psi\left((2N)^{-j}\left(\omega-\dfrac{\alpha}{2}\right)+\dfrac{\lambda-\sigma}{2}\right)}\\
&=\displaystyle\sum_{(j,\lambda)\ne(j,\sigma)\in I(\alpha)}\hat\psi\left(\dfrac{\omega}{(2N)^j}-(2N)^{-j}\dfrac{\alpha}{2}\right)
\overline{\hat\psi\left(\dfrac{\omega}{(2N)^j}-(2N)^{-j}\dfrac{\alpha}{2}+\dfrac{\lambda-\sigma}{2}\right)}\\
&=\displaystyle\sum_{(j,\lambda)\ne(j,\sigma)\in I(\alpha)}\hat\psi\left(\dfrac{\omega}{(2N)^j}-\dfrac{\lambda-\sigma}{2}\right)\overline{\hat\psi
\left(\dfrac{\omega}{(2N)^j}\right)}\\
&=\overline{\Delta_\alpha^-(\omega)}\qquad {\big(\text{By Eq. (2.13)}}.
\end{align*}

\parindent=0mm \vspace{.0in}
Therefore
\begin{align*}
\sum_{\alpha \in \Gamma\setminus\{0\}} \big|\Delta_\alpha^+(\omega)\big|=\sum_{\alpha \in \Gamma\setminus\{0\}} \big|\Delta_\alpha^{-}(\omega)\big|.\tag{2.15}
\end{align*}

\parindent=0mm \vspace{.0in}
Replacing $\xi+{\alpha}/{2}$ by $\omega$ in the last integration of (2.14), we derive from (2.14) and (2.15) that
\begin{align*}
\Big|R_\psi(f)\Big|&\le\left\{\displaystyle\sum_{\alpha \in \Gamma\setminus\{0\}}\int_{\mathbb R}\left|\hat f(\xi)\right|^2
\left|\Delta_\alpha^+(\xi)\right|d\xi\right\}^{1/2}\left\{\displaystyle
\sum_{\alpha \in \Gamma\setminus\{0\}}\int_{\mathbb R}\left|\hat f(\omega)\right|^2\left|\overline{\Delta_\alpha^-(\omega)}\right|d\omega\right\}^{1/2}\\
&=\displaystyle\int_{\mathbb R}\left|\hat f(\xi)\right|^2\left\{\displaystyle\sum_{\alpha \in \Gamma\setminus\{0\}}
\left|\Delta_\alpha^+(\xi)\right|\right\}d\xi.\tag{2.16}
\end{align*}

\parindent=0mm \vspace{.0in}
Therefore, by (2.16) and (2.4), we have
\begin{align*}
\int_{\mathbb R}\left|\hat f(\xi)\right|^2\left\{\Delta_0^+(\xi)-\sum_{\alpha \in \Gamma\setminus\{0\}}\left|\Delta_\alpha^+(\xi)\right|\right\}d\xi
 \le \sum_{j \in \mathbb Z}\sum_{\lambda \in \Lambda}\left|\big\langle f, \psi_{j, \lambda}\big\rangle\right|^2,\tag{2.17}
\end{align*}

\parindent=0mm \vspace{.0in}
and
\begin{align*}
\sum_{j \in \mathbb Z}\sum_{\lambda \in \Lambda}\left|\big\langle f, \psi_{j, \lambda}\big\rangle\right|^2
\le\int_{\mathbb R}\left|\hat f(\xi)\right|^2\left\{\Delta_0^+(\xi)+\sum_{\alpha \in \Gamma\setminus\{0\}}\left|\Delta_\alpha^+(\xi)\right|\right\}d\xi,
\end{align*}

\parindent=0mm \vspace{.0in}
or, equivalently
$$\sum_{j \in \mathbb Z}\sum_{\lambda \in \Lambda}\left|\big\langle f, \psi_{j, \lambda}\big\rangle\right|^2
\le \int_{\mathbb R}\left|\hat f(\xi)\right|^2\left\{\sum_{\alpha \in \Gamma}\left|\Delta_\alpha^+(\xi)\right|\right\}d\xi.\eqno (2.18)$$

\parindent=0mm \vspace{.1in}
Taking infimum in (2.17) and supremum in (2.18), respectively, we obtain again that

$$C_\psi\big\|f\big\|^2_2\le \sum_{j \in \mathbb Z}\sum_{\lambda \in \Lambda}\left|\big\langle f, \psi_{j, \lambda}\big\rangle\right|^2\le D_\psi\big\|f\big\|^2_2.$$

\parindent=0mm \vspace{.1in}
The proof of Theorem 2.4 is complete.\qquad\fbox

\parindent=0mm \vspace{.2in}

{\bf{ Remark 2.5.}} The frame bounds in Theorem 2.4 are better than that of Sharma and Manchanda [15, Theorem 3.1]. In fact, we have
the following:
\begin{align*}
A&=\displaystyle\inf_{\xi \in [1, 2N]}\left\{\sum_{j \in \mathbb Z}\left|\hat \psi\big((2N)^{-j}\xi\big)\right|^2-\sum_{j \in \mathbb Z}\sum_{\ell \ne 0}\left|\hat \psi\big((2N)^{-j}\xi\big)\overline{\hat \psi\big((2N)^{-j}\xi+\ell/2\big)}\right|\right\}\\
&=\displaystyle\inf_{\xi \in [1, 2N]}\left\{\sum_{j \in \mathbb Z}\left|\hat \psi\left(\dfrac{\xi}{(2N)^j}\right)\right|^2-\sum_{\alpha \in \Gamma\setminus\{0\}}\sum_{(j,\lambda)\ne(j,\sigma)\in I(\alpha)}\left|\hat \psi\left(\dfrac{\xi}{(2N)^j}\right)\overline{\hat\psi\left(\dfrac{\xi}{(2N)^j}+\dfrac{\lambda-\sigma}{2}\right)}\right|\right\}\\
&\le\displaystyle\inf_{\xi \in [1, 2N]}\left\{\sum_{j \in \mathbb Z}\left|\hat \psi\left(\dfrac{\xi}{(2N)^j}\right)\right|^2-\sum_{\alpha \in \Gamma\setminus\{0\}}\left|\sum_{(j,\lambda)\ne(j,\sigma)\in I(\alpha)}\hat \psi\left(\dfrac{\xi}{(2N)^j}\right)\overline{\hat\psi\left(\dfrac{\xi}{(2N)^j}+\dfrac{\lambda-\sigma}{2}\right)}\right|\right\}\\
&=\displaystyle\inf_{\xi\in[1, 2N]}\Big\{\Delta_0^+(\xi)-\sum_{\alpha \in \Gamma\setminus\{0\}}\left|\Delta_\alpha^+(\xi)\right|\Big\}\\
&=C_\psi,
\end{align*}

\parindent=0mm \vspace{.0in}
and
\begin{align*}
B&=\displaystyle\sup_{\xi \in [1, 2N]}\left\{\sum_{j \in \mathbb Z}\sum_{\ell \in \mathbb Z}\left|\hat \psi\big((2N)^{-j}\xi\big)\overline{\hat \psi\big((2N)^{-j}\xi+\ell/2\big)}\right|\right\}\\
&=\displaystyle\sup_{\xi \in [1, 2N]}\left\{\sum_{\alpha \in \Gamma\setminus\{0\}}\sum_{(j,\lambda)\ne(j,\sigma)\in I(\alpha)}\left|\hat \psi\left(\dfrac{\xi}{(2N)^j}\right)\overline{\hat\psi\left(\dfrac{\xi}{(2N)^j}+\dfrac{\lambda-\sigma}{2}\right)}\right|\right\}\qquad\qquad\qquad\\
&\ge\displaystyle\sup_{\xi \in [1, 2N]}\left\{\sum_{\alpha \in \Gamma\setminus\{0\}}\left|\sum_{(j,\lambda)\ne(j,\sigma)\in I(\alpha)}\hat \psi\left(\dfrac{\xi}{(2N)^j}\right)\overline{\hat\psi\left(\dfrac{\xi}{(2N)^j}+\dfrac{\lambda-\sigma}{2}\right)}\right|\right\}\\
&=\displaystyle\sup_{\xi\in[1, 2N]}\Big\{\sum_{\alpha \in \Gamma\setminus\{0\}}\left|\Delta_\alpha^+(\xi)\right|\Big\}\\
&=D_\psi.
\end{align*}

\parindent=0mm \vspace{.0in}
With the notations in (2.12) and (2.13), we define  new sets as
$$ \Pi_\alpha^+=\text{ess}\sup\Big\{\left|\Delta_\alpha^+(\xi)\right|: \xi \in [1, 2N]\Big\},\qquad \Pi_\alpha^{-}=\text{ess}\sup\Big\{\left|\Delta_\alpha^{-}(\xi)\right|:
 \xi \in [1, 2N]\Big\}.$$

\parindent=0mm \vspace{.1in}
{\bf{ Theorem 2.6.}} {\it Suppose $\psi \in L^2(\mathbb R)$ such that}
\begin{align*}
E_\psi&=\displaystyle ess\inf_{\xi \in [1, 2N]}\left\{\sum_{j \in \mathbb Z}\left|\hat \psi\left(\dfrac{\xi}{(2N)^j}\right)\right|^2\right\}-
\sum_{\alpha \in \Gamma\setminus\{0\}}\Big[\Pi_\alpha^+\,\Pi_\alpha^-\Big]^{1/2}>0\\
F_\psi&=\displaystyle ess\sup_{\xi \in [1, 2N]}\left\{\sum_{j \in \mathbb Z}\left|\hat \psi\left(\dfrac{\xi}{(2N)^j}\right)\right|^2\right\}+
\sum_{\alpha \in \Gamma\setminus\{0\}}\Big[\Pi_\alpha^+\,\Pi_\alpha^-\Big]^{1/2}< \infty.
\end{align*}

\parindent=0mm \vspace{.0in}
{\it Then $\big\{\psi_{j,\lambda}:j \in \mathbb Z,\lambda \in \Lambda\big\}$ is a wavelet frame for $L^2(\mathbb R)$ with bounds $E_\psi$ and $F_\psi$.}

\parindent=0mm \vspace{.2in}
{\it Proof.} By equation (2.14), we have
\begin{align*}
\Big|R_\psi(f)\Big|&\le\displaystyle\sum_{\alpha \in \Gamma\setminus\{0\}}\left\{\int_{\mathbb R}
\left|\hat f(\xi)\right|^2\left|\Delta_\alpha^+(\xi)\right|d\xi\right\}^{1/2}
\left\{\int_{\mathbb R}\left|\hat f\left(\xi+\dfrac{\alpha}{2}\right)\right|^2\left|\Delta_\alpha^+(\xi)\right|d\xi\right\}^{1/2}\\
&=\displaystyle\sum_{\alpha \in \Gamma\setminus\{0\}}\left\{\int_{\mathbb R}\left|\hat f(\xi)\right|^2\left|\Delta_\alpha^+
(\xi)\right|d\xi\right\}^{1/2}\left\{\int_{\mathbb R}
\left|\hat f(\omega)\right|^2\left|{\Delta_\alpha^-(\omega)}\right|d\xi\right\}^{1/2}\\
&\le\displaystyle\sum_{\alpha \in \Gamma\setminus\{0\}}\Big[\Pi_\alpha^+\,\Pi_\alpha^-\Big]^{1/2}\int_{\mathbb R}\left|\hat f(\xi)\right|^2 d\xi.
\end{align*}

\parindent=0mm \vspace{.0in}
Proceeding similarly as in Theorem 2.3, we have

$$\int_{\mathbb R}\left|\hat f(\xi)\right|^2\left\{\sum_{j \in \mathbb Z}\left|\hat \psi\left(\dfrac{\xi}{(2N)^j}\right)\right|^2-
\sum_{\alpha \in \Gamma\setminus\{0\}}\Big[\Pi_\alpha^+\,\Pi_\alpha^-\Big]^{1/2}\right\}d\xi
\le \sum_{j \in \mathbb Z}\sum_{\lambda \in \Lambda}\left|\big\langle f, \psi_{j, \lambda}\big\rangle\right|^2,$$

\parindent=0mm \vspace{.1in}
and
$$\sum_{j \in \mathbb Z}\sum_{\lambda \in \Lambda}\left|\big\langle f, \psi_{j, \lambda}\big\rangle\right|^2
\le \int_{\mathbb R}\left|\hat f(\xi)\right|^2\left\{\sum_{j \in \mathbb Z}\left|\hat \psi\left(\dfrac{\xi}{(2N)^j}\right)\right|^2+
\sum_{\alpha \in \Gamma\setminus\{0\}}\Big[\Pi_\alpha^+\,\Pi_\alpha^-\Big]^{1/2}\right\}d\xi.$$

\parindent=0mm \vspace{.1in}
These two inequalities imply that

$$E_\psi\big\|f\big\|^2_2\le \sum_{j \in \mathbb Z}\sum_{\lambda \in \Lambda}\left|\big\langle f, \psi_{j, \lambda}\big\rangle\right|^2\le F_\psi\big\|f\big\|^2_2.$$

\parindent=0mm \vspace{.1in}
This completes the proof of Theorem 2.6. \qquad\fbox

\parindent=0mm \vspace{.2in}


\parindent=0mm \vspace{.2in}

{\bf{References}}

\begin{enumerate}

{\small {

\item O. Christensen, {\it An Introduction to Frames and Riesz Bases}. Birkh\"{a}user, Boston, 2003.

\item C.K. Chui and X. Shi, Inequalities of Littlewood-Paley type for frames and wavelets, {\it SIAM J. Math. Anal.} 24 (1993), 263-277.

\item I. Daubechies, {\it Ten Lectures on Wavelets}, SIAM, Philadelphia, 1992.

\item I. Daubechies, A. Grossmann and Y. Meyer, Painless non-orthogonal expansions,  {\it J. Math. Phys.} 27(5) (1986), 1271-1283.

\item L. Debnath and F.A. Shah, {\it Wavelet Transforms and Their Applications}, Birkh\"{a}user, New York, 2015.

\item R.J. Duffin and A.C. Shaeffer, A class of nonharmonic Fourier series, {\it Trans. Amer. Math. Soc.} 72 (1952), 341-366.

\item B. Fuglede, Commuting self-adjoint partial different operators and a group theoretic problem,  {\it J. Funct. Anal.}  16  (1974), 101-121.

\item J.P. Gabardo and M.Z Nashed, Nonuniform multiresolution analysis and spectral pairs, {\it J. Funct. Anal.} 158  (1998), 209-241.

\item J.P. Gabardo and X. Yu, Wavelets associated with nonuniform multiresolution analysis and one-dimensional spectral pairs, {\it J. Math. Anal. Appl.}  323 (2006), 798-817.

\item D. Li, G. Wu and X. Yang, Unified conditions for wavelet frames,  {\it Georgian Math. J.} 18 (2011), 761-776.

\item F.A. Shah and Abdullah, Nonuniform multiresolution analysis on local fields of positive characteristic,  {\it Compl. Anal.  Opert. Theory.}
 9 (2015), 1589-1608.

\item F.A. Shah and M.Y. Bhat, Vector-valued nonuniform multiresolution analysis on local fields, {\it Int. J.  Wavelets, Multiresolut. Inf. Process.} 13(4) (2015), Article ID: 1550029.

\item F.A. Shah and M.Y. Bhat, Nonuniform wavelet packets on local fields of positive characteristic, {\it Filomat}. 31(6) (2017), 1491-1505.

\item F.A. Shah and L. Debnath, Dyadic wavelet frames on a half-line using the Walsh-Fourier transform, {\it Integ. Transf. Special Funct.} 22(7) (2011), 477-486.

\item V. Sharma and Manchanda, Nonuniform wavelet frames in $L^2(\mathbb R)$, {\it Asian-European J. Math.}  8 (2015), Article ID: 1550034.

   }}

\end{enumerate}

\end{document}